\documentclass[]{interact}

\usepackage{epstopdf}
\usepackage[caption=false]{subfig}

\usepackage[numbers,sort&compress]{natbib}
\bibpunct[, ]{[}{]}{,}{n}{,}{,}
\renewcommand\bibfont{\fontsize{10}{12}\selectfont}
\makeatletter
\def\NAT@def@citea{\def\@citea{\NAT@separator}}
\makeatother

\theoremstyle{plain}
\newtheorem{theorem}{Theorem}[section]
\newtheorem{lemma}[theorem]{Lemma}
\newtheorem{corollary}[theorem]{Corollary}

\theoremstyle{definition}

\theoremstyle{remark}

\usepackage{url} 
\newcommand{\diag}{\mathop{\mathrm{diag}}}
\renewcommand\vec{\mathbf}

\begin{document}

\articletype{Original Article}

\title{Bilinear matrix equation characterizes \\Laplacian and distance matrices of weighted trees}

\author{
\name{Mikhail Goubko\textsuperscript{a}\thanks{CONTACT Mikhail Goubko. Email: mgoubko@mail.ru} and Alexander Veremyev\textsuperscript{b}}
\affil{\textsuperscript{a}V.A. Trapeznikov Institute of Control Sciences of RAS, \\
65 Profsoyuznaya, Moscow, Russia; \\
\textsuperscript{b}University of Central Florida,
Orlando, FL, 32816-2368, USA}
}

\maketitle

\begin{abstract}
It is known from algebraic graph theory that if $L$ is the Laplacian matrix of some tree $G$ with a vertex degree sequence $\vec{d}=(d_1, ..., d_n)^\top$ and $D$ is its distance matrix, then $LD+2I=(2\cdot\vec{1}-\vec{d})\vec{1}^\top$, where $\vec{1}$ is an all-ones column vector. We prove that if this matrix identity holds for the Laplacian matrix of some graph $G$ with a degree sequence $\vec{d}$ and for some matrix $D$, then $G$ is essentially a tree, and $D$ is its distance matrix. This result immediately generalizes to weighted graphs. If the matrix $D$ is symmetric, the lower triangular part of this matrix identity is redundant and can be omitted. Therefore, the above bilinear matrix equation in $L$, $D$, and $\vec{d}$ characterizes trees in terms of their Laplacian and distance matrices. Applications to the extremal graph theory (especially, to topological index optimization and to optimal tree problems) and to road topology design are discussed.
\end{abstract}


\section{Introduction}
Let $I$ denote an identity matrix, $J$ be an all-ones matrix, and $\vec{1}$ be an all-ones column vector of convinient dimension, while $\diag \vec{a}$ standing for a diagonal matrix with elements of the column vector $\vec{a}$ on its main diagonal.

We study simple undirected graphs without loops. Let $V(G)$ and $E(G)$ denote the sets of vertices and edges of the graph $G$ respectively. A typical graph considered below has $n$ vertices indexed from 1 to $n$. Edges are bi-element sets $\{i,j\}\subseteq V(G)$ denoted $ij\in E(G)$ for short. The \textit{adjacency matrix} $X(G)=(x_{ij}(G))_{i,j=1}^n$ of the graph $G$ is a symmetric binary matrix with $x_{ij}(G)=1$ if $ij\in E(G)$ and $x_{ij}(G)=0$ otherwise. A \textit{degree sequence} of the graph $G$ is the column vector $\vec{d}(G):=X(G)\cdot\vec{1}$. 

Unless otherwise stated, \textit{weighted graphs} are considered, i.e., a positive weight $w_{ij}>0$ is assigned to every edge $ij\in E(G)$ of the graph $G$. The weight matrix $W(G)=(w_{ij}(G))_{i,j=1}^n$ of the graph $G$ is a non-negative symmetric matrix with $w_{ij}(G)=w_{ij}$ if $ij\in E(G)$ and $w_{ij}(G)=0$ otherwise. If edge weights are selected from some positive symmetric matrix $M=(\mu_{ij})_{i,j=1}^n$, then $W(G)=M\odot X(G)$, where $\odot$ is the elementwise matrix multiplication. Unweighted graphs can be considered a special case of weighted graphs with unitary weights. 

For a graph $G$ with edge weights $w_{ij}(G)$, $ij\in E(G)$, let $G^{-1}$ denote the graph $G$ with all weights replaced with their reciprocals: $w_{ij}\left(G^{-1}\right):=\frac{1}{w_{ij}(G)}$, $ij\in E(G)$. Obviously, $\left(G^{-1}\right)^{-1}=G$ and for unweighted graphs $G^{-1}=G$.

The Laplacian matrix of a weighted graph $G$ is $L(G):=\diag (W(G)\vec{1})-W(G)$. The distance $d_{ij}(G)$ between vertices $i$ and $j$ of a weighted graph $G$ is the minimum weight of the path that connects vertices $i$ and $j$ in graph $G$, where the weight of a path is just a sum of its edges' weights. The distance matrix $D(G)$ of the graph $G$ is an $n\times n$ matrix with zeros along its diagonal and with its $(i, j)$-entry being equal to $d_{ij}(G)$.

A connected graph with $n$ vertices and $n-1$ edges is called a \textit{tree}. Following the seminal works of Graham, Pollak, and Lov\'asz \cite{graham1971addressing,GRAHAM197860}, Laplacian and distance matrices of weighted and unweighted trees  are extensively studied in the literature. Let $L:=L(G)$ be the Laplacian matrix of some weighted tree $G$ with a degree sequence $\vec{d}=\vec{d}(G)$, and let $D:=D\left(G^{-1}\right)$ be the distance matrix of $G^{-1}$. Firstly introduced by Bapat, Kirkland, and Neumann (\cite{BAPAT2005193}, Lemma 4.1), the identity 
\begin{equation}\label{eq_LD}
LD+2I=(2\cdot\vec{1}-\vec{d})\vec{1}^\top,
\end{equation}
is an attractive result that brings together in a single expression the main quantities that characterize a graph in algebraic graph theory.

Being a simple consequence of the weighted version \cite{BAPAT2005193} of the famous Graham-Lov\'asz's formulae \cite{GRAHAM197860} for the inverse of the distance matrix of a tree,
\begin{equation}\label{eq_invD}
D^{-1}=\frac{(2\cdot\vec{1}-\vec{d})(2\cdot\vec{1}-\vec{d})^\top}{2\sum_{ij\in E(G)}w_{ij}\left(G^{-1}\right)}-\frac{1}{2}L,
\end{equation}
the identity (\ref{eq_LD}) and its immediate corollary $LDL=-2L$ are widely used to establish algebraic properties of trees \cite{Balaji2007,BALAJI2007108,bapat2016squared}. 

In the next section we prove that the converse proposition is also valid, i.e., if identity (\ref{eq_LD}) holds for the Laplacian matrix of some graph $G$ and for some matrix~$D$, then $G$ is essentially a tree with the degree sequence $\vec{d}$ and $D$ is the distance matrix of the graph $G^{-1}$. It is also shown that in case of a symmetrix matrix $D$, equations in the lower-triangular part of the matrix equation (\ref{eq_LD}) are redundant and can be omitted.

In Section \ref{sec_apps} applications to the extremal graph theory and to topology design problems are outlined. The discussion on the possible extensions of the above results concludes.

\section{Results}

\begin{theorem}\label{th_LD1}
Let $L$ be the Laplacian matrix of a (weighted) graph $G$ with vertex degree sequence $\vec{d}=(d_1,...,d_n)^\top$. If for some real matrix $D$ the identity \emph{(\ref{eq_LD})} holds, then $G$ is a tree. Moreover, if $D$ has the zero diagonal, then $D$ is the distance matrix of the tree $G^{-1}$.
\begin{proof}
To prove that $G$ is a tree, let us multiply both sides of the equation (\ref{eq_LD}) by $\vec{1}^\top$ from the left and by $\vec{1}$ from the right, obtaining
\begin{equation}\label{eq_deg1}
\vec{1}^\top LD\vec{1} + 2n = (2n-\vec{1}^\top\vec{d})n,
\end{equation}
Since $L$ is a graph Laplacian matrix, it has an all-ones eigenvector $\vec{1}$ that corresponds to the zero eigenvalue, and, therefore, $\vec{1}^\top L\equiv 0$. Hence, from (\ref{eq_deg1}) it follows that $\vec{1}^\top\vec{d}=2(n-1)$, i.e., 
\begin{equation}\label{eq_tree_seq}
\sum_{i=1}^n d_i=2(n-1),
\end{equation} 
and $\vec{d}$ must be the vertex degree sequence of some tree.

Therefore, if graph $G$ is connected, then $G$ is a tree. By contradiction, assume that $G$ is disconnected. Then $L$ has an eigenvector $\vec{u}$ for the zero eigenvalue, such that $\vec{u}^\top L=0$ and $\vec{1}^\top \vec{u} =0$.
Multiplying the equation (\ref{eq_LD}) by $\vec{u}^\top$ from the left, we obtain $2\vec{u}^\top + \vec{u}^\top\vec{d}\vec{1}^\top = 0$. It follows that $u_i=u_j=-\vec{u}^\top\vec{d}/2$, which is impossible since $\vec{1}^\top u =0$. So, $G$ is a tree.

Let us prove that the matrix $D$ is equal to the distance matrix $D\left(G^{-1}\right)$. From \cite{BAPAT2005193,bapat2010graphs} it is known that if $G$ is a tree, the equation (\ref{eq_LD}) holds for $D = D\left(G^{-1}\right)$. Therefore,  
$L\Delta=0$, where $\Delta = (\Delta_{ij})_{i,j=1}^n := D-D\left(G^{-1}\right)$.

Multiplying $L\Delta=0$ from the left by the generalized inverse \cite{Gutman2004GeneralizedIO} $L^\dag$ of the Laplacian matrix $L$, and taking into account that $L^\dag L=I-\frac{1}{n}J$ \cite{bapat2004resistance,Balaji2007} (remember that $J$ is an all-ones matrix), we obtain that
$$\Delta=\frac{1}{n}J\Delta,$$
i.e., $\Delta_{ij}=\frac{1}{n}\sum_{k=1}^n \Delta_{ik}$. This means that $\Delta_{ij}=\Delta_{ii}$ for all $i,j=1,...,n$. Since $D$ has the zero diagonal, $\Delta_{ii}\equiv 0$. Therefore, $\Delta_{ij}=0$ and $D \equiv D\left(G^{-1}\right)$, which completes the proof.
\end{proof}
\end{theorem}

For unweighted graphs $G \equiv G^{-1}$, so, the following corollary holds.
\begin{corollary}\label{cor_LD1}
If $L$ is the Laplacian matrix of an unweighted graph $G$ with vertex degree sequence $\vec{d}=(d_1,...,d_n)^\top$, and for some real matrix $D$ 
the identity \emph{(\ref{eq_LD})} holds, then $G$ is a tree. If $D$ has the zero diagonal, then $D$ is a distance matrix of the tree $G$.
\end{corollary}

A symmetric $n\times n$ matrix with zero diagonal and positive off-diagonal elements is called a \textit{Eucidean distance matrix} (EDM) if $n$ points can be selected in some $p$-dimensional Euclidean space so that element $(i,j)$ of the matrix is equal to the squared Euclidean distance between the $i$-th and the $j$-th points. If, in addition, points can be selected on a $p$-dimensional hypersphere, the matrix is called a \textit{spherical} (or circulant) EDM.

For the next theorem we will need the following auxiliary result.
\begin{lemma}\label{lemma_CEDM}
The distance matrix $D$ of a (weighted) tree $G$ is a spherical EDM. 
\begin{proof}
Trees with positive edge weights are a special case of matrix-weighted trees \cite{BALAJI2007108}. Balaji and Bapat \cite{BALAJI2007108} proved that the distance matrix of a matrix-weighted tree is EDM. Bapat, Kirkland, and Neumann showed that $\det D\neq 0$, so $D$ is invertible \cite{BAPAT2005193}. Gower \cite{gower1985properties} showed that EDM $D$ is spherical if $\vec{1}^\top D^{-1}\vec{1} > 0$. Substituting the explicit expression (\ref{eq_invD}) for $D^{-1}$ and taking into account that $L\vec{1}\equiv 0$, we see that 
$$\vec{1}^\top D^{-1}\vec{1}=\frac{\left(\vec{1}^\top(2\cdot\vec{1}-\vec{d})\right)^2}{2\sum_{ij\in E(G)}w_{ij}\left(G\right)} =\frac{2}{\sum_{ij\in E(G)}w_{ij}\left(G\right)}>0.$$\end{proof}\end{lemma}

For a square matrix $A=(a_{ij})_{i,j=1}^n$ let $A^\urcorner=(a_{ij}^\urcorner)_{i,j=1}^n$ denote its upper triangular part (excluding the diagonal), i.e., $a_{ij}^\urcorner = a_{ij}$ if $i<j$, and $a_{ij}^\urcorner = 0$ otherwise.

\begin{theorem}\label{th_LD2}
If $L:=L(G)$ is a Laplacian matrix of some (weighted) graph $G$ with vertex degree sequence $\vec{d}=(d_1,...,d_n)^\top$, and for some real symmetric zero-diagonal matrix $D$, such that $D\ge W\left(G^{-1}\right)$, the identity 
\begin{equation}\label{eq_LDruc}
\left[LD + 2I - (2\cdot \vec{1}-\vec{d})\vec{1}^\top\right]^\urcorner = 0
\end{equation}
holds, then $G$ is a tree and the matrix $D=D\left(G^{-1}\right)$ is the distance matrix of $G^{-1}$.
\begin{proof}
The matrix equation (\ref{eq_LDruc}) is equivalent to the equation

\begin{equation}\label{eq_LDZ}
LD + 2I = (2\cdot \vec{1}-\vec{d})\vec{1}^\top + Z,
\end{equation}
where $Z$ is some lower-triangular matrix (which means that $Z^\urcorner \equiv 0$).

Let us introduce the spectral decomposition of the graph Laplacian matrix $L=U\Lambda_0 U^\top$. Here $\Lambda_0 := \diag{\vec{\lambda}_0}$, where $\vec{\lambda}_0 := (\lambda_1, ..., \lambda_{n-c}, 0, ..., 0)$ is the vector of eigenvalues enumerated in the descending order, and $U:=(\vec{u_1}, ..., \vec{u_n})$ is the orthogonal matrix whose columns $\vec{u_1}, ..., \vec{u_n}$ are the corresponding normalized eigenvectors. The Laplacian matrix of a graph $G$ has non-negative eigenvalues with $\lambda_n=0$, and the multiplicity $c$ of zero eigenvalue is equal to the number of connected components in graph $G$. 

Let us assume that graph $G$ is disconnected (i.e., $c\ge 2$) and prove that in this case the matrix equation (\ref{eq_LDZ}) has no solution. The Laplacian eigenvectors can always be selected so that: 
\begin{enumerate} 
\item $\vec{u}_n=\frac{1}{\sqrt{n}}(1,..., 1)^\top$, 
\item If $c\ge 2$ then 
$\vec{u}_{n-1}=\alpha \vec{1}(C_1)+\beta \vec{1}(C_2)$, 
where $C_1, C_2 \subset V(G)$ are disjoint vertex sets of the graph $G$ such that $C_1\cup C_2 = V(G)$, and $\vec{1}(C_1), \vec{1}(C_2)$ are their corresponding indicatory vectors, while $\alpha$ and $\beta$ are constants. Let us denote $n_1:=|C_1|>0$, $n_2 := |C_2|=n-n_1>0$.
\end{enumerate} 

Eigenvectors are orthonormal, so $\vec{u}_{n-1}^\top\vec{u}_n=0$, $\vec{u}_n^\top\vec{u}_n=1$, which implies that 
\begin{equation}\label{eq_alpha}
\alpha=\sqrt{\frac{n_2}{n_1n}},  \beta = -\sqrt{\frac{n_1}{n_2n}}.
\end{equation}

Let us multiply both sides of (\ref{eq_LDZ}) by the row-vector $\vec{u}_n^\top=\frac{1}{\sqrt{n}}\vec{1}^\top$ from the left. By the definition of a zero-eigenvalue eigenvector, $\vec{u}_n^\top L=0$, so, equality (\ref{eq_LDZ}) gives 
\begin{equation}\label{eq_1TZ}
\vec{1}^\top Z = \left(\sum_{i=1}^n d_i-2(n-1)\right)\vec{1}^\top.
\end{equation}
Introducing the shorthand notation $\vec{z}_i:=\vec{1}(C_i)^\top Z,i=1,2$, we can write an equivalent equation:
\begin{equation} \label{eq_LDZ1}
\vec{z}_1+\vec{z}_2 = \left(\sum_{i=1}^n d_i-2(n-1)\right)\vec{1}^\top. 
\end{equation}

Then, multiplying both sides of the equation (\ref{eq_LDZ}) by the row-vector $\vec{u}_{n-1}^\top$ from the left and taking into account orthogonality of eigenvectors, we have another system of equations: 
\begin{equation} \label{eq_LDZ2}
\alpha \vec{z}_1 + \beta \vec{z}_2 = 2\alpha\vec{1}(C_1)^\top + 2\beta\vec{1}(C_2)^\top + ( \alpha\sum_{i\in C_1}d_i + \beta\sum_{i\in C_2}d_i)\vec{1}^\top,
\end{equation}
which, after simplification, gives
\begin{equation} \label{eq_LDZ21}
\vec{z}_1 - \frac{n_1}{n_2} \vec{z}_2 = 2\vec{1}(C_1)^\top - 2\frac{n_1}{n_2} \vec{1}(C_2)^\top +\left ( \sum_{i\in C_1}d_i - \frac{n_1}{n_2} \sum_{i\in C_2}d_i\right)\vec{1}^\top.
\end{equation}

Combining (\ref{eq_LDZ1}) and (\ref{eq_LDZ21}), we evaluate
\begin{equation} \label{eq_z1}
\vec{z}_1 = \left(\sum_{i\in C_1}d_i-2n_1\right)\vec{1}^\top+2\cdot\vec{1}(C_1)^\top. 
\end{equation}

Without loss of generality, let us assume that $n\in C_2$. Since, by definition, $Z$ is a lower-diagonal matrix, it follows that $n$-th component of $\vec{z}_1$ is zero. Then, from (\ref{eq_z1}) it is clear that $\sum_{i\in C_1}d_i=2n_1$ and 
\begin{equation} \label{eq_z12}
\vec{z}_1 =2\cdot\vec{1}(C_1)^\top.
\end{equation} 

Let $m := \max\{i: i\in C_1\}$ be the maximum index of a vertex in $C_1$. From (\ref{eq_z12}) we conclude that $m$-th diagonal element of the matrix $Z$ is equal to $2$. So, writing down the $m$-th diagonal element of the matrix equation (\ref{eq_LDZ}), we have $(LD)_{mm}+2=2-d_m+2$, or, using notation $D=(d_{ij})_{i,j=1}^n$,
\begin{equation}\label{eq_LDmm}
\sum_{k=1}^n w_{mk}(G) d_{mm}-\sum_{k=1}^n w_{mk}(G)d_{km}=2-d_m.
\end{equation}

By the hypothesis, $d_{mm}\equiv0$, so (\ref{eq_LDmm}) is equivalent to 
$$\sum_{k=1}^n w_{mk}(G)d_{km}=d_m-2.$$
On the other hand, by the hypothesis, $D\ge W\left(G^{-1}\right)$, so $d_{km}\ge \frac{1}{w_{km}(G)}$ for all $km\in E(G)$. Since $w_{km} = w_{mk}$ and $w_{mk}(G)\equiv0$ for all $mk\notin E(G)$, 
$$d_m-2=\sum_{k=1}^n w_{mk}(G)d_{km}=\sum_{km\in E(G)} w_{km}(G)d_{km} \ge\sum_{km\in E(G)} \frac{w_{km}(G)}{w_{km}(G)}=d_m,$$ 
which makes a contradiction.

Therefore, if the matrix equation (\ref{eq_LDruc}) holds, then the graph $G$ is connected. 

By the Euler theorem, 
\begin{equation}\label{eq_euler}
\sum_{i=1}^n d_i \ge 2(n-1) 
\end{equation}
for any connected graph $G$. 

Since $Z=(z_{ij})_{i,j=1}^n$ is lower-triangular, from (\ref{eq_1TZ})  it follows that 
\begin{equation}\label{eq_znn1}
z_{nn}=\sum_{i=1}^n d_i-2(n-1). 
\end{equation}
On the other hand, $n$-th diagonal element of the matrix equation (\ref{eq_LDZ}) gives
\begin{equation}\label{eq_LDnn}
\sum_{k=1}^n w_{nk}(G) d_{nn}-\sum_{k=1}^n w_{nk}(G)d_{kn}+2=2-d_n+z_{nn}.
\end{equation}
As before, taking into account that $d_{nn}\equiv0$ and $d_{kn}\ge \frac{1}{w_{kn}(G)}$ for $kn\in E(G)$, we conclude that $\sum_{k=1}^n w_{nk}(G)d_{kn}\ge d_n$ and
$$z_{nn}=d_n-\sum_{k=1}^n w_{nk}(G)d_{kn}\le 0.$$

Substituting (\ref{eq_znn1}), we obtain that $\sum_{i=1}^n d_i\le 2(n-1)$, which, in combination with (\ref{eq_euler}), gives $\sum_{i=1}^n d_i = 2(n-1)$. Therefore, by the Euler theorem, $G$ has $n-1$ edges. It is already shown that $G$ is connected, so $G$ is a tree.

Let us prove that the matrix $D$ is equal to the distance matrix $D\left(G^{-1}\right)$ of the graph $G^{-1}$. From \cite{BAPAT2005193,bapat2010graphs} it is known that if $G$ is a tree, the equation (\ref{eq_LD}) holds for $D = D\left(G^{-1}\right)$. Subtracting the equation (\ref{eq_LD}) for $D = D\left(G^{-1}\right)$ from the equation (\ref{eq_LDZ}) we obtain  
the equation $L(D-D\left(G^{-1}\right))=Z$.

Balaji and Bapat \cite{Balaji2007} showed that the generalized inverse $L^\dag$ of the Laplacian matrix $L$ of the tree $G$ is evaluated as $L^\dag=\left(I-\frac{1}{n}J\right)D\left(G^{-1}\right)\left(I-\frac{1}{n}J\right)$. 
Multiplying both sides of the equation $L(D-D\left(G^{-1}\right))=Z$ by $L^\dag$ from the left, we obtain
\begin{equation}\label{eq_LdagL}
L^\dag L(D-D\left(G^{-1}\right))=\left(I-\frac{1}{n}J\right)D\left(G^{-1}\right)\left(I-\frac{1}{n}J\right)Z.
\end{equation}
Since $G$ is a tree, from (\ref{eq_1TZ}) we know that $\vec{1}^\top Z=0$ and, consequently, $JZ=0$. Taking into account that $L^\dag L=I-\frac{1}{n}J$, the equation (\ref{eq_LdagL}) reduces to
\begin{equation}\label{eq_final}
\left(I-\frac{1}{n}J\right)\left(D-D\left(G^{-1}\right)+\frac{1}{2}D\left(G^{-1}\right)Z\right)=0,
\end{equation}
which is equivalent to
\begin{equation}\label{eq_D}
\Delta+\frac{1}{2}D\left(G^{-1}\right)Z = \vec{1}\cdot \vec{a}^\top,
\end{equation}
where $\Delta=(\Delta_{ij})_{i,j=1}^n:=D-D\left(G^{-1}\right)$ and $\vec{a}=(a_i)_{i=1}^n$ is an unknown column vector.

Since $Z=(z_{ij})_{i,j=1}^n$ is a lower-triangular matrix, $z_{in}=0$ for all $i=1,...,n-1$. Also, the matrix equation $\vec{1}^\top Z = 0$ implies that $z_{nn}$ is also equal to zero. Therefore, the $n$-th column of matrix $D\left(G^{-1}\right)Z$ is null and, writing the $n$-th column of the matrix equation (\ref{eq_D}), we have $\Delta_{ij}=d_{in}-d_{in}\left(G^{-1}\right)=a_n$ for all $i=1,...,n$. 

By definition, $d_{nn}=d_{nn}\left(G^{-1}\right)=0$, so $a_n=0$. Thus, we conclude that $d_{in}=d_{in}\left(G^{-1}\right)$ and, by symmetry, $d_{ni}=d_{ni}\left(G^{-1}\right)$ for all $i=1,...,n$.

Recalling that $z_{k,n-1}\equiv 0$ for all $k\le n-2$, $(n-1)$-th column of the equation (\ref{eq_D}) gives 
\begin{equation}\label{eq_n_1_system}\begin{cases}
\Delta_{1,n-1}+\frac{1}{2}\left[d_{1,n-1}\left(G^{-1}\right)z_{n-1,n-1} + d_{1n}\left(G^{-1}\right)z_{n,n-1}\right]=a_{n-1}\\
...\\
\Delta_{n-2,n-1}+\frac{1}{2}\left[d_{n-2,n-1}\left(G^{-1}\right)z_{n-1,n-1} + d_{n-2,n}\left(G^{-1}\right)z_{n,n-1}\right]=a_{n-1}\\
\frac{1}{2}d_{n-1,n}\left(G^{-1}\right)z_{n,n-1}=a_{n-1}\\
\frac{1}{2}d_{n,n-1}\left(G^{-1}\right)z_{n-1,n-1}=a_{n-1}
\end{cases}
\end{equation} 

From the matrix equation $\vec{1}^\top Z = 0$ we obtain the equation $z_{n-1,n-1}+z_{n,n-1}=0$. Taking into account that distances are positive, so $d_{n-1,n}\left(G^{-1}\right)>0$, and, by symmetry, $d_{n-1,n}\left(G^{-1}\right)=d_{n,n-1}\left(G^{-1}\right)$, the last two equations of the system (\ref{eq_n_1_system}) are written as
$$d_{n-1,n}\left(G^{-1}\right)z_{n,n-1}=2a_{n-1}, -d_{n-1,n}\left(G^{-1}\right)z_{n,n-1}=2a_{n-1}.$$
These equations immediately give $a_{n-1}=0$, and, consequently, $z_{n,n-1}=-z_{n-1,n-1}=0$. 
Substituting to the system (\ref{eq_n_1_system}), we obtain $\Delta_{i,n-1}=0$ and $d_{i,n-1}=d_{i,n-1}\left(G^{-1}\right)$ for all $i=1,...,n$. 

A similar procedure is repeated recursively for columns $n-2, n-3, ..., 1$ of the equation (\ref{eq_D}). For instance, the $j$-th column of (\ref{eq_D}) can be written as
\begin{equation}\label{eq_j_system_full}
\vec{\Delta}_j + \frac{1}{2}D\left(G^{-1}\right)\vec{z}_j=a_j\cdot \vec{1},
\end{equation}
where $\vec{\Delta}_j:=(\Delta_{ij})_{i=1}^{n}$ and $\vec{z}_j:=(z_{ij})_{i=1}^{n}$.

Using expression (\ref{eq_invD}) for the inverse of the distance matrix, 
we write 
\begin{equation}\label{eq_z_j_full}
\vec{z}_j=\left(\frac{(2\cdot\vec{1}-\vec{d})(2\cdot\vec{1}-\vec{d})^\top}{\sum_{ij\in E(G)}w_{ij}\left(G^{-1}\right)}-L\right)(a_j\cdot \vec{1}-\vec{\Delta}_j).
\end{equation}

Multiplying both sides of (\ref{eq_z_j_full}) by $\vec{1}^\top$ from the left, and taking into account that $\vec{1}^\top \vec{z}_j = 0$, $\vec{1}^\top L=0$, an expression for $a_j$ follows after simplifications:
\begin{equation}\label{eq_a_j}
a_j =\frac{1}{2}(2\cdot\vec{1}-\vec{d})^\top\vec{\Delta}_j.
\end{equation}

On the other hand, when arriving to the $j$-th column, we already know that $\Delta_{ij}=0$ for all $i> j$. Since $d_{jj}=d_{jj}\left(G^{-1}\right)=0$, we also know that $\Delta_{jj}=0$. As the matrix $Z$ is lower-triangular, $z_{ij}\equiv 0$ for all $i< j$. With the notation $\vec{\Delta}_j':=(\Delta_{ij})_{i=1}^{j-1}$ and $\vec{z}_j':=(z_{ij})_{i=j}^{n}$ for the unknown variables, equations (\ref{eq_j_system_full}),(\ref{eq_a_j}) reduce to the following system of linear equations:
\begin{equation}\label{eq_j_system}
\begin{cases}
\vec{\Delta}_j' + \frac{1}{2}D_{12}\vec{z}_j'=a_j\cdot \vec{1},\\
\frac{1}{2}D_{22}\vec{z}_j'=a_j\vec{1},\\
a_j =\frac{1}{2}(2\cdot\vec{1}-\vec{d}_1)^\top\vec{\Delta}_j',
\end{cases}
\end{equation}
where $D\left(G^{-1}\right) = \begin{pmatrix}
D_{11} & D_{12}\\
D_{21} & D_{22}
\end{pmatrix}$, $D_{11}$ is $(j-1)\times(j-1)$ matrix, $\vec{d}_1:=(d_i)_{i=1}^{j-1}$, and all-ones column vectors have compatible dimension.

It is a well-known fact that $D_{22}$ non-singular, as it is a principal submatrix of order at least
2 of the distance matrix of a weighted tree \cite{BAPAT2005193}. Hence, from the second line of the system (\ref{eq_j_system}), $\vec{z}_j'=2a_jD_{22}^{-1}\vec{1}$. 

Multiplying both sides of this equation by $\vec{1}^\top$ from the left and taking into account that $\vec{1}^\top \vec{z}_j'=0$, we have 
\begin{equation}\label{eq_a_or_D} 
0=2a_j\vec{1}^\top D_{22}^{-1}\vec{1}. 
\end{equation}

Lemma \ref{lemma_CEDM} says that $D\left(G^{-1}\right)$ is a spherical EDM. As a principal submatrix of a spherical EDM, the matrix $D_{22}$ is also a spherical EDM with the same assignment of points in the same Euclidean space. 

Gower \cite{gower1985properties} showed that $\vec{1}^\top D_{22}^{-1}\vec{1}\neq 0$ if $D_{22}$ is a non-singular spherical EDM. So, from (\ref{eq_a_or_D}) it follows that $a_j=0$. Substituting $a_j=0$ into the system (\ref{eq_j_system}) immediately gives the desired equalities $\vec{z}_j=0$ and $\Vec{\Delta}_j=0$ (and, therefore, $d_{ij}=d_{ij}\left(G^{-1}\right)$ for all $i=1,....,n$).
\end{proof}
\end{theorem}

\begin{corollary}\label{cor_LD2}
If $L$ is the Laplacian matrix of an unweighted graph $G$ with vertex degree sequence $\vec{d}=(d_1,...,d_n)^\top$, and for some real symmetric zero-diagonal matrix  $D\ge J-I$ the identity \emph{(\ref{eq_LDruc})} holds, then $G$ is a tree and $D$ is its distance matrix.
\end{corollary}

\section{Applications}\label{sec_apps}
According to Theorem \ref{th_LD1}, $n^2$ equations of the matrix equation (\ref{eq_LD}) simultaneously assure that: 
\begin{enumerate}
\item the Laplacian matrix $L$ has rank $n-1$, so the graph $G$ is connected, 
\item the degree sequence $\vec{d}$ sums up to $2(n-1)$, so $G$ is a tree, 
\item the matrix $D$ conforms $L$, so that $D$ is the distance matrix of $G^{-1}$.
\end{enumerate}

Distance matrices of trees are a special case of EDM \cite{Balaji2007}. So, the result of Theorem \ref{th_LD1} is in line with the well-known one-to-one correspondence between EDMs and Laplacian matrices (i.e., positive semidefinite matrices of rank $n-1$). However, Theorem \ref{th_LD1} does not impose neither the computationally intractable rank constraints on $L$ nor any constraint on $D$ to assure that it is a EDM. The equations of (\ref{eq_LD}) become attractive equality constraints for optimization problems that concern trees.

In particular, (\ref{eq_LD}) is linear both in the elements of $L$ and $D$ matrices. If $L$ is fixed, (\ref{eq_LD}) becomes a system of linear equations on $D$, which has no solution unless $L$ is a Laplacian matrix of some (weighted) tree $G$, in which case its solution is the distance matrix of $G^{-1}$. If $D$ and $\vec{d}$ are fixed, (\ref{eq_LD}) is a linear system on $L$, which is inconsistent unless $D$ is a distance matrix of a tree $G$ with the degree sequence $\vec{d}$, in which case its solution is the Laplacian matrix of $G^{-1}$.

It is not the most efficient way to calculate neither the distance nor the Laplacian matrix of a tree. However, if elements of both $L$ and $D$ matrices are free variables, the expression (\ref{eq_LD}) gives a system of bilinear equations that characterize (weighted) trees in terms of their Laplacian and distance matrices.

To be more specific, let us consider a collection of graphs with edge weights selected from some positive real symmetric matrix $M = (\mu_{ij})_{i,j=1}^n$. Every such graph is characterized by $\frac{n(n-1)}{2}$ binary variables $x_{ij}$ for $1\le i<j\le n$ that determine its adjacency matrix $X=(x_{ij})_{i,j=1}^n$. Its weight matrix $W:=M\odot X$, the degree sequence $\vec{d}:=X\cdot\vec{1}$, and the Laplacian matrix $L:=\diag(W\vec{1})-W$ are linear expressions in $X$. Let us also introduce  $\frac{n(n-1)}{2}$ real variables $d_{ij}$, $1\le i<j\le n$, that define a symmetric zero-diagonal matrix $D=(d_{ij})_{i,j=1}^n$. 

Then, by Theorem \ref{th_LD1}, any solution of the system (\ref{eq_LD}) of $n^2$ bilinear equations in $\frac{n(n-1)}{2}$ binary variables $x_{ij}$ and $\frac{n(n-1)}{2}$ real variables $d_{ij}$, $1\le i<j\le n$, corresponds to some tree $G$ with adjacency matrix $x_{ij}(G)=x_{ij}$, reciprocal edge weights $w_{ij}(G)=\frac{1}{\mu_{ij}}$, $ij\in E(G)$, and the distance matrx $D(G)=D$. And vice versa, Bapat, Kirkland, and Neumann \cite{BAPAT2005193} showed that any tree $G$ with edge weights $w_{ij}(G)=\frac{1}{\mu_{ij}}$, $ij\in E(G)$, corresponds to $\frac{n(n-1)}{2}$ binaries $x_{ij}=x_{ij}(G)$ and $\frac{n(n-1)}{2}$ real numbers $d_{ij}=d_{ij}(G)$, $1\le i<j\le n$, which determine matrices $X$ and $D$ that reduce $n^2$ equations (\ref{eq_LD}) to the identities.

Therefore, any extremal tree problem with a criterion being linear in $X$ and $D$ can be written as a mixed-integer bilinear program (MIBP). Although equality constraints in (\ref{eq_LD}) involving binary variables are non-convex, they can be converted into linear constraints (see \cite{GurobiMIBP}) by introducing one auxiliary real variable and adding $4$ linear inequality constraints for each of $\frac{n(n-1)(n-2)}{6}$ independent bilinear terms $x_{ik}d_{kj}$, $1\le i < k< j\le n$. The problem thus becomes a mixed-integer linear program (MILP), which is efficienttly handled by the modern optimization software like Gurobi. 

Moreover, Theorem \ref{th_LD2} says that $\frac{n(n+1)}{2}$ equations in the lower-diagonal part of the matrix equation (\ref{eq_LD}) can be replaced by $\frac{n(n-1)}{2}$ simpler inequality constraints of the form $d_{ij}\ge \frac{x_{ij}}{\mu_{ij}}$, $1\le i<j\le n$.

The set of admissible trees can be futher limited by additional constraints employing $X$ and $D$. For example, we can limit ourselves to the trees with a given number of pendent vertices or a given degree sequence by adding linear constraints on the components of the vector $\vec{d} = X\cdot\vec{1}$. An edge $ij$ can be banned or enforced with conditions $x_{ij}=0$ or $x_{ij}=1$ respectively. Inequality constraints on $d_{ij}$ allow to impose an upper limit on the diameter of considered trees or on the eccentricity of their vertices \cite{Dankelmann2020BoundingT}.

The problems of tree topology design that can be reduced to a MILP as described above include (but not limited to) optimization of degree- \cite{goubko2014note,goubko2014degree} and distance-based topological indices \cite{xu2014survey} (also known as molecular descriptors in mathematical chemistry \cite{todeschini2008handbook}), especially, the Wiener index \cite{SUN2019438}, its version for vertex-weighted trees \cite{goubko2016minimizing,goubko2018maximizing}, and the weighted average distance of a tree \cite{goubko2020lower}. Some of these problems are still unsolved, and MILP solutions can be a good starting point of the analysis.

The above-listed examples of topological indices involve weighted sums of distances for unweighted graphs. But Theorems \ref{th_LD1} and \ref{th_LD2} allow for the more general, and still understudied, case of the sum of weighted distances \cite{cai2019sum}.  The total graph weight can also be a criterion for the MILP, as in the classical \textit{minimum spanning tree} (MST) problem and its NP-hard versions (e.g., \textit{degree-constrained} MST, hop-constrained MST \cite{GOUVEIA1996178,PIRKUL2003126,GOUVEIA2020364} or the Steiner tree problems)

Probably the simplest cost function for road network topology design (RNTD) problems is a combination of the total graph weight $C\odot X$ (representing the road investment costs for the construction cost matrix $C$) and of the weighted sum $M\odot D$ of weighted distances (modelling the total traveling cost in a graph with weighted edges given an origin-destination matrix $M$) \cite{jia2019review}. Thus, the problem of tree-shaped RNTD, as explained above, also reduces to a MILP.

\section{Conclusion}
In this article, we have shown that the system of matrix equations (\ref{eq_LD}) can be used to characterize the Laplacian and the distance matrices of weighted trees. For trees with $n$ vertices there are $\frac{n(n-1)}{2}$ binary and $\frac{n(n-1)}{2}$ real independent variables in the system (\ref{eq_LD}) of $n^2$ bilinear equality constraints or, alternatively, in the system (\ref{eq_LDruc}) of $\frac{n(n-1)}{2}$ bilinear equality and $\frac{n(n-1)}{2}$ linear inequality constraints that can be used in MIBP settings of tree topology design problems. These MIBPs reduce to MILPs by introducing $O(n^3)$ auxiliary variables and adding $O(n^3)$ linear inequality constraints.

Among the mixed integer programs that give both graph topology and the distance matrix, this setting is among the most compact in terms of the number of variables and constraints. Alternative approaches require $O(n^3)$ (binary) variables with $O(n^4)$ linear constraints \cite{GOUVEIA1996178}, $O(n^4)$ binaries with $O(n^4)$ constraints \cite{balakrishnan1992using}, or $O(n^5)$ binaries with $O(n^5)$ constraints \cite{PIRKUL2003126} under different flow-based approaches, and $O(n^4)$ binaries with $O(n^4)$ constraints  under the recent path-based approach \cite{veremyev2015critical,mukherjee2017minimum,diaz2019robust}.

The results of this article can be improved in several directions:
\begin{enumerate}
\item It is an open question whether inequalities $D\ge W(G^{-1})$ in the conditions of Theorem \ref{th_LD2} are essential or can be relaxed. 
\item Balaji and Bapat introduced graphs with matrix weights \cite{Balaji2007} and showed that their distance matrices inherit many properties of those for weighted graphs. So we might conjecture that both Teorem \ref{th_LD1} and Theorem \ref{th_LD2} directly generalize to the trees with matrix weights.
\item Unicyclic graphs are a natural generalization of trees. As shown by Bapat, Kirkland, and Neumann \cite{BAPAT2005193}, unicyclic graphs have interesting algebraic characteristics, and probably might be characterized similarly to trees.
\item Attractive algebraic properties of shortest-path distances in trees immediately generalize to \textit{resistance distances} \cite{klein1993resistance} in general weighted graphs with cycles. In particular, Bapat proved \cite{bapat2004resistance} analogs of both the expression (\ref{eq_invD}) for the inverse distance matrix and of the matrix equation (\ref{eq_LD}) for arbitrary connected weighted graphs (a similar theory for unweighted graphs was earlier developed by Xiao and Gutman \cite{xiao2003resistance}). Therefore, one can pose a question wherther a compact characterization exists of general connected graphs similar to the one for trees introduced in this article, and whether this characterization can be used in extremal graph problems for resistance-distance-based topological indices.
\end{enumerate}
These generalizations can be the subject of the future work.

\section*{Funding}
This work was supported by the Russian Foundation for Basic Research (RFBR) [18-07-01240].

\bibliographystyle{tfnlm}
\bibliography{references}

\begin{thebibliography}{10}
\providecommand{\url}[1]{\normalfont{#1}}
\providecommand{\urlprefix}{Available from: }

\bibitem{graham1971addressing}
Graham~RL, Pollak~HO. On the addressing problem for loop switching. The Bell
  System Technical Journal. 1971;\hspace{0pt}50(8):2495--2519.

\bibitem{GRAHAM197860}
Graham~RL, Lov\'asz~L. Distance matrix polynomials of trees. Advances in
  Mathematics. 1978;\hspace{0pt}29(1):60 -- 88.
  \urlprefix\url{http://www.sciencedirect.com/science/article/pii/0001870878900051}.

\bibitem{BAPAT2005193}
Bapat~RB, Kirkland~SJ, Neumann~M. On distance matrices and {L}aplacians. Linear
  Algebra and its Applications. 2005;\hspace{0pt}401:193 -- 209. Special Issue
  in honor of Graciano de Oliveira;
  \urlprefix\url{http://www.sciencedirect.com/science/article/pii/S0024379504002599}.

\bibitem{Balaji2007}
Balaji~R, Bapat~RB. Block distance matrices. ELA The Electronic Journal of
  Linear Algebra [electronic only]. 2007;\hspace{0pt}16:435--443.
  \urlprefix\url{http://eudml.org/doc/129304}.

\bibitem{BALAJI2007108}
Balaji~R, Bapat~RB. On {E}uclidean distance matrices. Linear Algebra and its
  Applications. 2007;\hspace{0pt}424(1):108 -- 117. Special Issue in honor of
  Roger Horn;
  \urlprefix\url{http://www.sciencedirect.com/science/article/pii/S0024379506002631}.

\bibitem{bapat2016squared}
Bapat~RB, Sivasubramanian~S. Squared distance matrix of a tree: inverse and
  inertia. Linear Algebra and its Applications. 2016;\hspace{0pt}491:328--342.

\bibitem{bapat2010graphs}
Bapat~RB. Graphs and matrices. Vol.~27. Springer; 2010.

\bibitem{Gutman2004GeneralizedIO}
Gutman~I, Xiao~W. Generalized inverse of the {L}aplacian matrix and some
  applications. Bulletin: Classe Des Sciences Mathematiques Et Natturalles.
  2004;\hspace{0pt}129:15--23.

\bibitem{bapat2004resistance}
Bapat~RB. Resistance matrix of a weighted graph. MATCH Communications in
  Mathematical and in Computer Chemistry. 2004;\hspace{0pt}50:73--82.

\bibitem{gower1985properties}
Gower~JC. Properties of {E}uclidean and non-{E}uclidean distance matrices.
  Linear Algebra and its Applications. 1985;\hspace{0pt}67:81--97.

\bibitem{GurobiMIBP}
Achterberg~T. Products of variables in mixed integer programming
  [\url{https://www.gurobi.com/wp-content/uploads/2019/07/2019-07-23_Products_of_Variables_in_Mixed_Integer_Programming_KA.pdf}];
  2019. Accessed: 2020-08-01.

\bibitem{Dankelmann2020BoundingT}
Dankelmann~P, Dossou-Olory~AAV. Bounding the $k$-{S}teiner {W}iener and
  {W}iener-type indices of trees in terms of eccentric sequence. arXiv:
  Combinatorics. 2020;\hspace{0pt}.

\bibitem{goubko2014note}
Goubko~M, R{\'e}ti~T. Note on minimizing degree-based topological indices of
  trees with given number of pendent vertices. MATCH Commun Math Comput Chem.
  2014;\hspace{0pt}72(3):633--639.

\bibitem{goubko2014degree}
Goubko~M, Gutman~I. Degree-based topological indices: {O}ptimal trees with
  given number of pendents. Applied Mathematics and Computation.
  2014;\hspace{0pt}240:387--398.

\bibitem{xu2014survey}
Xu~K, Liu~M, Das~KC, et~al. A survey on graphs extremal with respect to
  distance-based topological indices. MATCH Commun Math Comput Chem.
  2014;\hspace{0pt}71(3):461--508.

\bibitem{todeschini2008handbook}
Todeschini~R, Consonni~V. Handbook of molecular descriptors. Vol.~11. John
  Wiley \& Sons; 2008.

\bibitem{SUN2019438}
Sun~Q, Ikica~B, \v{S}krekovski~R, et~al. Graphs with a given diameter that
  maximise the {W}iener index. Applied Mathematics and Computation.
  2019;\hspace{0pt}356:438 -- 448.
  \urlprefix\url{http://www.sciencedirect.com/science/article/pii/S0096300319302267}.

\bibitem{goubko2016minimizing}
Goubko~M. Minimizing {W}iener index for vertex-weighted trees with given weight
  and degree sequences. Match: communications in mathematical and computer
  chemistry. 2016;\hspace{0pt}75(1):3--27.

\bibitem{goubko2018maximizing}
Goubko~M. Maximizing {W}iener index for trees with given vertex weight and
  degree sequences. Applied Mathematics and Computation.
  2018;\hspace{0pt}316:102--114.

\bibitem{goubko2020lower}
Goubko~M, Kuznetsov~A. Lower bound for the cost of connecting tree with given
  vertex degree sequence. Journal of Complex Networks.
  2020;\hspace{0pt}8(2):cnz031.

\bibitem{cai2019sum}
Cai~Q, Li~T, Shi~Y, et~al. Sum of weighted distances in trees. Discrete Applied
  Mathematics. 2019;\hspace{0pt}257:67--84.

\bibitem{GOUVEIA1996178}
Gouveia~L. Multicommodity flow models for spanning trees with hop constraints.
  European Journal of Operational Research. 1996;\hspace{0pt}95(1):178 -- 190.
  \urlprefix\url{http://www.sciencedirect.com/science/article/pii/0377221795000909}.

\bibitem{PIRKUL2003126}
Pirkul~H, Soni~S. New formulations and solution procedures for the hop
  constrained network design problem. European Journal of Operational Research.
  2003;\hspace{0pt}148(1):126 -- 140.
  \urlprefix\url{http://www.sciencedirect.com/science/article/pii/S0377221702003661}.

\bibitem{GOUVEIA2020364}
Gouveia~L, Leitner~M, Ljubi\'c~I. A polyhedral study of the diameter
  constrained minimum spanning tree problem. Discrete Applied Mathematics.
  2020;\hspace{0pt}285:364 -- 379.
  \urlprefix\url{http://www.sciencedirect.com/science/article/pii/S0166218X20302614}.

\bibitem{jia2019review}
Jia~GL, Ma~RG, Hu~ZH. Review of urban transportation network design problems
  based on {C}ite{S}pace. Mathematical Problems in Engineering.
  2019;\hspace{0pt}2019.

\bibitem{balakrishnan1992using}
Balakrishnan~A, Altinkemer~K. Using a hop-constrained model to generate
  alternative communication network design. ORSA Journal on Computing.
  1992;\hspace{0pt}4(2):192--205.

\bibitem{veremyev2015critical}
Veremyev~A, Prokopyev~OA, Pasiliao~EL. Critical nodes for distance-based
  connectivity and related problems in graphs. Networks.
  2015;\hspace{0pt}66(3):170--195.

\bibitem{mukherjee2017minimum}
Mukherjee~T, Veremyev~A, Kumar~P, et~al. The minimum edge compact spanner
  network design problem. arXiv preprint arXiv:171204010. 2017;\hspace{0pt}.

\bibitem{diaz2019robust}
Diaz~C, Nikolaev~A, Perla~A, et~al. Robust communication network formation: a
  decentralized approach. Computational Social Networks.
  2019;\hspace{0pt}6(1):1--30.

\bibitem{klein1993resistance}
Klein~DJ, Randi{\'c}~M. Resistance distance. Journal of mathematical chemistry.
  1993;\hspace{0pt}12(1):81--95.

\bibitem{xiao2003resistance}
Xiao~W, Gutman~I. {R}esistance distance and {L}aplacian spectrum. Theoretical
  Chemistry Accounts. 2003;\hspace{0pt}110(4):284--289.

\end{thebibliography}

\end{document}